\documentclass[12pt,reqno,a4paper]{amsart}
\usepackage{amssymb,amscd}

\begin{document}

\let\kappa=\varkappa
\let\eps=\varepsilon
\let\phi=\varphi

\def\Z{\mathbb Z}
\def\R{\mathbb R}
\def\C{\mathbb C}
\def\Q{\mathbb Q}

\def\OO{\mathcal O}
\def\CP{\C{\mathrm P}}
\def\RP{\R{\mathrm P}}
\def\conj{\overline}
\def\Beta{\mathrm{B}}
\def\p{\partial}
\def\dist{\mathop{\mathrm{dist}}}
\def\supp{\mathop{\mathrm{supp}}}

\renewcommand{\Im}{\mathop{\mathrm{Im}}\nolimits}
\renewcommand{\Re}{\mathop{\mathrm{Re}}\nolimits}
\newcommand{\codim}{\mathop{\mathrm{codim}}\nolimits}
\newcommand{\id}{\mathop{\mathrm{id}}\nolimits}
\newcommand{\Aut}{\mathop{\mathrm{Aut}}\nolimits}

\newtheorem*{thm}{Theorem}
\newtheorem*{prop}{Proposition}

\theoremstyle{definition}
\newtheorem*{rems}{Remarks}
\newtheorem*{exm}{Example}

\renewcommand{\thesubsection}{\arabic{subsection}}

\title{Finite unions of  balls in $\C^n$ are rationally convex}

\author{Stefan Nemirovski}
\address{%
Steklov Mathematical Institute;\hfill\break
\strut\hspace{8 true pt} Ruhr-Universit\"at Bochum}
\email{stefan@mi.ras.ru}
\thanks{The author was supported by the Russian Foundation for Basic
Research (grant no.\ 05-01-00981), the programme ``Leading
Scientific Schools of Russia'' (grant no.\ 9429.2006.1), the
programme ``Contemporary Problems of Theoretical Mathematics'' of
the Russian Academy of Sciences, and the project SFB/TR 12 of the
Deutsche Vorschungsgemeinschaft.}

\maketitle

\begin{flushright}
\small
\begin{tabular}{l}
-- Is this a prologue, or the posy of a ring?\\
-- 'Tis brief, my lord.
\end{tabular}\\[2pt]
{\it Hamlet}, III:2
\end{flushright}

\bigskip
\noindent
{\bf 1.} A compact set $K\subset\C^n$ is called {\it polynomially convex\/}
if for any point $z\notin K$ there exists a polynomial $P$ such that
$|P(z)|>\max\limits_{\xi\in K}|P(\xi)|$. Replacing polynomials by rational functions,
one gets the definition of a {\it rationally convex\/} compact set. These notions
are interesting, in particular, because any function holomorphic in a neighbourhood
of a polynomially (respectively, rationally) convex set can be uniformly on this set
approximated by polynomials (respectively, rational functions).

An old problem asks whether any finite union of  disjoint closed balls in $\C^n$ is polynomially convex.
It is known only that the answer is positive for at most three balls~\cite{Ka}. In this note,
we show that the {\it rational\/} convexity of any such union follows almost immediately
from the results of Julien Duval and Nessim Sibony~\cite{DS}.

\begin{thm}
Any union of finitely many disjoint closed balls in $\C^n$ is rationally convex.
\end{thm}

Note that it follows from the construction of the examples in \cite{Ka} and \cite{Kh}
and the argument principle that this statement is false for polydiscs and complex ellipsoids in $\C^3$.

\medskip
\noindent
{\bf 2.}
Let us first recall that according to Theorem~1.1 in~\cite{DS}, {\it if $\omega$ is
a non-negative $d$-closed $(1,1)$-form on $\C^n$ such that $\C^n\setminus\supp\omega$ is
relatively compact in $\C^n$, then for any $s>0$, the set
$\{z\in\C^n\mid \dist(z,\supp\omega)\ge s\}$ is rationally convex}.
We shall need the following corollary of this result.

\begin{prop}
Let $\phi$ be a strictly plurisubharmonic function on an open subset $U\subset\C^n$
such that its Levi form $dd^c\phi$ extends to a positive $d$-closed $(1,1)$-form on the whole $\C^n$.
If the set $K_\phi=\{z\in U\mid \phi(z)\le 0\}$ is compact, then  it is rationally convex.
\end{prop}

\begin{proof}
Fix a small $\eps>0$
and consider a smooth convex non-decreasing function  $f:\R\to\R$  such that
$f'(t)\equiv 0$ for $t\le\eps$, $f'(t)>0$ for $\eps<t<2\eps$, and $f'(t)\equiv 1$ for $t\ge 2\eps$.
Note that
$$
dd^c f(\phi)=f'(\phi) dd^c\phi + f''(\phi)d\phi\wedge d^c\phi\ge f'(\phi)dd^c\phi.
$$
Hence, if we set
$$
\omega_\eps=
\begin{cases}
dd^c\phi & \text{on } \C^n\setminus \{\phi\le 2\eps\} ,\\
dd^c f(\phi)&  \text{on } \{\phi\le 2\eps\},
\end{cases}
$$
then $\omega_\eps$ satisfies the conditions of the Duval--Sibony theorem
and its support $\supp\omega_\eps$ is precisely $\C^n\setminus\{\phi <\eps\}$.
As $\eps$ and $s=s(\eps)$ can be chosen arbitrarily close to zero, we see that $K_\phi$ is the intersection of
rationally convex sets and hence rationally convex itself.
\end{proof}

\begin{rems}
1$^\circ$ It follows from other results in \cite{DS} that any rationally convex compact set in~$\C^n$ has a fundamental system of neighbourhoods of the form $\{\phi<0\}$, where $\phi$ satisfies
the assumptions of the proposition.

\smallskip
\noindent
2$^\circ$ For comparison, note that a compact set $K\subset\C^n$ is polynomially convex
if and only if it has a fundamental system of neighbourhoods of the form $\{\phi<0\}$,
where now $\phi$ is an exhausting strictly plurisubharmonic function on  the whole $\C^n$.
\end{rems}

\smallskip
\noindent
{\bf 3.}
We can now prove the theorem. Let $\conj B(a_j,r_j)=\{z\in\C^n\mid\|z-a_j\|^2\le r_j^2\}$, $j=1,\dots,N$,
be a collection of pairwise disjoint closed balls. In a neighbourhood of their union $\bigsqcup \conj B(a_j,r_j)$,
consider the function $\phi$ that is equal to
$$
\phi_j(z)=\|z-a_j\|^2-r_j^2
$$
near each $\conj B(a_j,r_j)$.
Then $K_\phi=\bigsqcup \conj B(a_j,r_j)$ and
$$
dd^c\phi=\frac{i}{2}\sum\limits_{k=1}^n dz_k\wedge d\conj z_k
$$
is the standard flat K\"ahler form on $\C^n$. Hence, $\bigsqcup \conj B(a_j,r_j)$
is rationally convex by the proposition from \S 2.\qed

\end{document}